\documentclass{article}
\usepackage{graphicx} 

\title{On a property of the polynomials \(x^n+(1-x)^n+a^n\)}
\author{Hayk Karapetyan\\Yerevan State University\\1 Alek Manukyan, Yerevan, Armenia}
\date{}

\usepackage{amsmath, amssymb, amsthm, cancel}

\renewcommand{\Re}{\text{Re }}
\newtheorem{definition}{Definition}
\newtheorem{theorem}{Theorem}
\newtheorem{lemma}{Lemma}
\newtheorem{proposition}{Proposition}
\newtheorem{corollary}{Corollary}
\newtheorem{conjecture}{Conjecture}
\numberwithin{equation}{section}
\numberwithin{theorem}{section}
\numberwithin{definition}{section}
\numberwithin{lemma}{section}
\numberwithin{corollary}{section}

\newcommand{\ZZ}{\mathbb{Z}}

\usepackage{tikz}

\begin{document}
	
	\maketitle
	
	\begin{abstract}
		The polynomials \(x^n + (1-x)^n + a^n \) arise naturally from FLT (Fermat's Last Theorem). We formulate a conjecture about them which is a generalization of FLT. We investigate the complex roots of these polynomials, and our main result is that in the cases \(|a|\leq \frac12\) and \(a=-1\), they lie on an explicitly given curve while 'filling in' that curve. We hypothesize that this property can be generalized to hold for other \(a\) as well.
	\end{abstract}
	
	\section{Introduction}
	
	We investigate the polynomials \(K_{a,n}:=x^n+(1-x)^n+a^n\), where \(a\neq 0\) is a rational parameter. These polynomials arise naturally from the following reformulation of FLT:
	\begin{proposition}
		\(K_{a,n}\) has a rational root for some rational \(a\notin\{0, -1\}\) and odd \(n\) iff \(X^n+Y^n=Z^n,~n>2\) has a non-trivial positive integer solution.
	\end{proposition}
	\begin{proof}
		\((\Rightarrow)\) Let \(x\) be the root. If \(x=0\) or \(1-x=0\), then \(1+a^n=0\) which cannot be the case as \(a\neq -1\). Hence, \(x,1-x,-a\) are rational and nonzero. They satisfy \(x^n + (1-x)^n = (-a)^n\). FLT can be equivalently formulated as the nonexistence of \textit{rational} solutions to \(X^n+Y^n=Z^n\) with \(XYZ\neq 0\) (see \cite{Edwards1996}, p. 3 for proof), hence the proposition.
		
		\((\Leftarrow)\) Observe that if a non-trivial solution to FLT exists for some integer \(n\), then a solution exists for all divisors of \(n\). \(n\) cannot equal to \(2^m\) as it would imply a solution for \(n=4\), something that was ruled out by Fermat himself (see \cite{Edwards1996}, pp. 9-10). This means that there is a solution for an \(n\) odd prime. In that case,  
		\[X^n+Y^n=Z^n \implies \left( \frac{X}{X+Y}\right)^n + \left(\frac{Y}{X+Y}\right)^n + \left(-\frac{Z}{X+Y}\right)^n = 0,\]
		so \(K_{a,n}\) has a rational root \(x = \frac{X}{X+Y}\) for \(a=-\frac{Z}{X+Y}\). Evidently \(a\neq 0\) and \(a\neq -1\) as \((X+Y)^n > X^n + Y^n = Z^n\).
		
	\end{proof}
	
	FLT just says that \(K_{a,n}\) doesn't have rational roots for odd \(n\) and \(a\in\mathbb{Q} \setminus\{0, -1\}\). One can consider a more general question: are these polynomials even reducible over \(\mathbb{Q}\)?
	
	\begin{conjecture}
		The polynomials \(K_{a,n}\), where \(a\) is a rational different from \(0\) and \(-1\), are irreducible over \(\mathbb{Q}\).
	\end{conjecture}
	
	We have performed computer calculations which confirm this conjecture for \(a\) and \(n\) such that \(a=\pm\frac{p}{q},a \neq -1,~ 0<p,q<200\) and \(n<100\). We will analyze the location of the roots of \(K_{a,n}\) on the complex plane as it generally is a helpful tool in analyzing the number-theoretic properties of a polynomial. The results we obtain work for any \(n\).
	
	\begin{definition}
		We will say that a polynomial set \(P_n, n\in \mathbb{N}\) \textbf{localizes} on a union of curves \(\gamma\) if the roots of \(P_n\) lie on \(\gamma\) and form an everywhere dense set on \(\gamma\).
	\end{definition}
	
	Fix \(a\in\{-1\}\cup \left[ -\frac12, \frac12\right] \). Our main result is that \(K_{a,n}\) localizes on an explicitly given union of curves: for \(|a|\leq \frac12\), it is the line \(\Re x = \frac12\), and for \(a=-1\), it is the dark blue curve in Fig. 1. Other examples of localizing polynomial sets include \(x^n-1\) (on the unit circle) and Chebyshev polynomials of the first and second kind \(T_n(x)\) and \(U_n(x)\) (on the segment \([-1,1]\))
	
	Investigating the reducibility of \(K_{a,n}\) when \(a=-1\) is of independent interest. It turns out that \(K_{-1,n}\) may have the factors \(x,x-1,\) or \(x^2-x+1\), with multiplicities that will be exactly given depending on \(n\). If we remove these "trivial" factors from the canonical decomposition of \(K_{-1,n}\), it is unclear if the remaining polynomial is irreducible or not. Computer calculations show that it is indeed irreducible for \(n<4000\).
	
	\section{The case \(|a|\leq\frac12\)}
	
	\begin{theorem}\label{maintrick}
		Fix \(a\in\mathbb{R}\). Then at least \(\left\lfloor\frac n2\right\rfloor-\left\lceil \frac{n}{\pi}\arccos\min\left( 1, \frac{1}{2|a|}\right) \right\rceil\) many roots of \(K_{a,n}\) lie on the line \(\Re x = \frac12\) not strictly above the point \(A\), where \(A\) is the upper point on the line with modulus \(\max \left( \frac12,|a|\right) \). Those roots form an everywhere dense set on that curve.
	\end{theorem} 
	\begin{proof}
		Consider the variable written in the form $x=\frac{1}{2} + \frac{1}{2}i\tan\theta$, where $\theta\in D:= \left[ \arccos\min\left( 1,\frac1{2|a|}\right) ,\frac\pi2\right) $. Note that for the lowest value of \(\theta\), 
		
		\begin{multline*}
			|x|^2 = x\bar{x} = \left(\frac{1}{2} + \frac{1}{2}i\tan\theta\right)\left(\frac{1}{2} - \frac{1}{2}i\tan\theta\right)=\frac1{4\cos^2\theta} =\\= \frac1{4\min\left(1,\dfrac{1}{2|a|}\right)^2 }  = \max\left(\frac{1}{4}, |a|\right)^2,
		\end{multline*}
		
		so the map \(\theta \mapsto \frac{1}{2} + \frac{1}{2}i\tan\theta \) indeed maps \(D\) to the ray \([A,A + i\infty)\). Now,
		
		\begin{multline*}
			K_{a,n}(x) = \left( \frac{1}{2} + \frac{1}{2}i\tan\theta\right) ^n + \left(\frac{1}{2} - \frac{1}{2}i\tan\theta\right) ^n + a^n =\frac{(\cos\theta+i\sin\theta)^n}{(2\cos\theta)^n}+\\+\frac{(\cos\theta-i\sin\theta)^n}{(2\cos\theta)^n}+a^n = \frac{2\cos n\theta}{(2\cos\theta)^n}+a^n=\frac{2\cos n\theta + (2a\cos\theta)^n}{(2\cos\theta)^n}.
		\end{multline*}
		
		It is sufficient to find the zeroes of \(f_n(\theta):=2\cos n\theta + (2a\cos\theta)^n\). Consider \(f_n\) defined on \(\bar{D}\) (the topological closure). Observe that since either \(|a|\leq \frac12\) or \(\theta \geq \arccos \frac1{2|a|}\), $|2a\cos \theta|\leq 1$. Hence, \(f_n(\theta)\) has the same sign as $\cos n\theta$ when $\cos n\theta=\pm 1$, equivalently \(\theta = \frac{k\pi}{n},~k\in\ZZ\). The values of \(k\) for which \(\theta\in \bar{D}\) are the integers in \(\left[\frac{n}{\pi}\arccos\min\left( 1, \frac{1}{2|a|}\right), \frac{n}{2}\right]\). There are \(\left\lfloor\frac n2\right\rfloor-\left\lceil \frac{n}{\pi}\arccos\min\left( 1, \frac{1}{2|a|}\right) \right\rceil+1\) possible values for such \(k\). \(\cos n\theta\) has different signs for successive points \(\theta = \frac{k\pi}{n}\) and \(\theta = \frac{(k+1)\pi}{n}\), so \(f_n\) has different signs as well. Since \(f_n\) is continuous and real-valued, there are zeros between these successive values of \(\theta\) (these zeroes are all in \(D\) as \(k=\frac{n}2\) does not yield a zero). Therefore, there are at least \(\left\lfloor\frac n2\right\rfloor-\left\lceil \frac{n}{\pi}\arccos\min\left( 1, \frac{1}{2|a|}\right) \right\rceil\) zeroes of \(f_n\) (and hence \(K_{a,n}\)). 
		
		For any interval \(\left[ \frac{k\pi}{n}, \frac{(k+1)\pi}{n}\right)\subset D \), there exists a root of \(f_n\), so all the roots form an everywhere dense set in \(D\). This everywhere dense set is mapped to an everywhere dense set on the ray \([A, A+i\infty)\) by the homeomorphism \(\theta \mapsto \frac{1}{2} + \frac{1}{2}i\tan\theta \).
		
	\end{proof}
	\begin{corollary}
		If \(|a|\leq\frac12\), \(K_{a,n}\) localizes on the line \(\Re x = \frac12\).
	\end{corollary}
	\begin{proof}
		Fix \(a\) and \(n\). According to the theorem, there are at least \(\left\lfloor\frac n2\right\rfloor \) roots in the upper half-plane. Since \(K_{a,n}\) has real coefficients, the conjugates of its roots are also roots. Thus, we obtain at least \(2\left\lfloor\frac n2\right\rfloor\) roots on the line \(\Re x = \frac12\)\footnote{The conjugates are distinct from the originals as \(x=\frac12\) is not a root: it corresponds to \(\theta = 0\), and \(f_n(\theta) = 2 + (2a)^n \neq 0\).}. But \(2\left\lfloor\frac n2\right\rfloor=\deg K_{a,n}\), so there are no other roots.
		
		When changing \(n\), the set of roots is everywhere dense above \(\frac12\). New roots obtained by conjugation form an everywhere dense set below \(\frac12\). Hence, the set of roots is everywhere dense on the whole line.
	\end{proof}
	
	\section{The case \(a=-1\)}
	
	Denote \(K_n(x) = K_{-1,n}(x)\). We will first consider the multiplicities of the complex zeros of $K_n$.
	
	Note that $K_n'=n(x^{n-1}-(1-x)^{n-1})$, and let \(G := \gcd(K_n, K_n')\). 
	\begin{multline}\label{eq:modtrick}
		x^{n-1} \equiv (1-x)^{n-1} \pmod G\implies 0\equiv x^n + (1-x)^{n-1}(1-x) + (-1)^n \equiv\\ \equiv x^{n}+x^{n-1}(1-x) + (-1)^n = x^{n-1} +(-1)^n\pmod G,
	\end{multline}
	
	so \(G\mid x^{n-1} +(-1)^n\).
	In the complex plane, the roots of the latter lie on the unit circle. The roots of \(K_n'\) satisfy \(|x|^{n-1}=|1-x|^{n-1}\Rightarrow |x|=|1-x|\), so they lie on the line \(\Re x = \frac12\). The roots of \(G\) must lie on both of these curves. They intersect at $\omega$ and $\bar\omega$, where $\omega$ is the primitive sixth root of $1$. Thus $G = (x-\omega)^k(x-\bar\omega)^{k_1}$, and since \(G\) has real coefficients, \(G = (x-\omega)^k(x-\bar\omega)^{k}\).
	To find out the multiplicities of those roots in \(K_n\), note that
	\begin{gather*}
		K_n(\omega) := \omega^n + \bar\omega^n + (-1)^n = \begin{cases}
			3, &\text{ if } n\equiv 0\pmod 6,\\
			0, &\text{ if } n\equiv 1\pmod 6,\\
			0, &\text{ if } n\equiv 2\pmod 6,\\
			-3, &\text{ if } n\equiv 3\pmod 6,\\
			0, &\text{ if } n\equiv 4\pmod 6,\\
			0, &\text{ if } n\equiv 5\pmod 6
		\end{cases},\\
		K_n'(\omega) = (n-1)\left( \omega^{n-1} - \bar\omega^{n-1}\right)  = \begin{cases}
			0, &\text{ if } n\equiv 1\pmod 6,\\
			i\sqrt3(n-1), &\text{ if } n\equiv 2\pmod 6,\\
			0, &\text{ if } n\equiv 4\pmod 6,\\
			-i\sqrt3(n-1), &\text{ if } n\equiv 5\pmod 6
		\end{cases},\\
		K_n''(\omega) = (n-1)(n-2)\left( \omega^{n-2} + \bar\omega^{n-2}\right) = \begin{cases}
			(n-1)(n-2), &\text{ if } n\equiv 1\pmod 6,\\
			-(n-1)(n-2), &\text{ if } n\equiv 4\pmod 6
		\end{cases}.
	\end{gather*}
	Thus, if $n\equiv0\pmod3$, \(\omega\) is not a root of \(K_n\), if \(n\equiv 2\pmod 3\), it is a simple root, and if $n\equiv1\pmod3$, it is a double root. 
	
	Note as well that $0$ and $1$ are roots of $K_n$ if \(n\) is odd. They are not roots of \(K'_n\), so they are simple. To talk about roots other than \(0,1,\omega,\bar\omega\), we will denote \(\tilde{K}_n\) the polynomial obtained by removing the factors \(x,x-1,x-\omega,x-\bar\omega\) from the canonical decomposition of \(K_n\). Denote \(d_n = \deg \tilde{K}_n\). Considering the fact that \(\deg K_n = n\) if \(n\) is even and \(n-1\) if \(n\) is odd, it is easy to calculate that \[
	d_n = \begin{cases}
		n-7,\text{ if } n\equiv 1\pmod 6\\
		n-(n\mod6), \text{ otherwise}
	\end{cases}.
	\] 
	Hence, \(d_n\) is divisible by \(6\).
	
	Call \(L\) the union of two rays  \((\omega, \omega+i\infty)\cup(\bar{\omega}, \bar\omega-i\infty)\). Call \(A_1\) the circular arc from \(\omega\) to \(\bar{\omega}\) passing through \(0\) (the center of this circle is at \(1\)). Call \(A_2\) the arc from \(\omega\) to \(\bar{\omega}\) passing through \(1\) (see Fig. 1).
	
	\begin{theorem}
		\(K_n\) localizes on \(L\cup A_1 \cup A_2\).
	\end{theorem}
	\begin{proof}
		It is equivalent to \(\tilde{K}_n\) localizing on \(L\cup A_1 \cup A_2\) as removing the roots \(0,1,\omega,\bar\omega\) has no effect on the location of other roots and everywhere density.
		
		Theorem \ref{maintrick} implies (with taking conjugates of the roots) that there are at least \(2\left( \left\lfloor\frac{n}{2}\right\rfloor - \left\lceil\frac{n}{3}\right\rceil\right) \) many roots of $\tilde{K}_n$ on \(L\). Straightforward checking shows \( \left\lfloor\frac{n}{2}\right\rfloor - \left\lceil\frac{n}{3}\right\rceil = \frac{d_n}{6}\). Now observe that 
		\[
		K_n(x)=K_n(1-x)=(-x)^{n}K_n\left( \frac1{x}\right). 
		\]
		This implies that the roots of \(\tilde{K}_n(x)\) are mapped to other roots under the maps $x\longmapsto 1-x$ and $x\longmapsto \frac1x$. Since the map \(x\longmapsto \frac1x\) is a geometric inversion with center \(0\) and radius \(1\) followed by a reflection across the real axis, it maps the \(\frac{d_n}{3}\) roots on \(L\) to \(A_1\)\footnote{See \cite{inversion} pp. 77-83 for the basic theory of geometric inversion. The inversion of the line \(\Re x=\frac12\) is a circle passing through zero. It should also pass through \(\omega\) and \(\bar\omega\) since those remain fixed.}. Then the map \(x\longmapsto 1-x\) (a central symmetry across \(\frac{1}{2}\)) maps those to \(A_2\). In total, we have found \(d_n\) many roots on the curves, so there can be no more.
		
		Since the maps $x\longmapsto 1-x$ and $x\longmapsto \frac1x$ are homeomorphisms on \(A_1\) and \(L\) respectively, they map the set of everywhere dense roots on \(L\) to a set of everywhere dense roots on \(A_2\) and \(A_1\).
		
		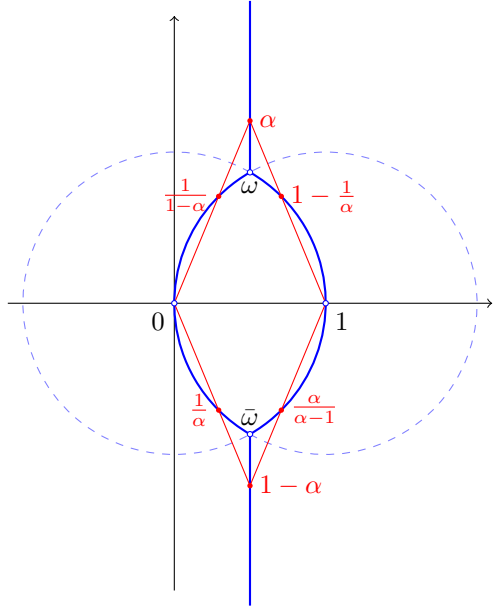
\begin{figure}[h]
			\centering
			\begin{tikzpicture}[scale = 2]
				\draw[->] (-1.1,0) -- (2.1,0);
				\draw[->] (0,-1.9) -- (0,1.9);
				
				\coordinate[label=225:$0$] (0) at (0,0);
				\coordinate[label=360-45:$1$] (1) at (1,0);
				\coordinate[label=below:$\omega$] (om) at (0.5,0.866);
				\coordinate[label=above:$\bar\omega$] (bom) at (0.5,-0.866) {};
				
				\draw[thick,blue] (bom) arc (-60:60:1);
				\draw[thick,blue] (om) arc (120:240:1);
				\draw[dashed, blue!50!white] (om) arc (60:300:1);
				\draw[dashed, blue!50!white] (bom) arc (-120:120:1);
				
				\draw[blue, thick] (om) -- (0.5, 2);
				\draw[blue, thick] (bom) -- (0.5, -2);
				
				\fill[red] (0.707, 0.707) circle (0.5pt) node[anchor=west] {$1-\frac1\alpha$};
				\fill[red] (0.293, 0.707) circle (0.5pt) node[anchor=east] {$\frac1{1-\alpha}$};
				\fill[red] (0.707, -0.707) circle (0.5pt) node[anchor=west] {$\frac\alpha{\alpha-1}$};
				\fill[red] (0.293, -0.707) circle (0.5pt) node[anchor=east] {$\frac{1}{\alpha}$};
				\fill[red] (0.5, 1.207) circle (0.5pt) node[anchor=west] {$\alpha$};
				\fill[red] (0.5, -1.207) circle (0.5pt) node[anchor=west] {$1-\alpha$};
				\draw[red] (1) -- (0.5, 1.207) -- (0) -- (0.5, -1.207) -- cycle;
				
				\filldraw[white, draw=blue] (0) circle (0.5pt);
				\filldraw[white, draw=blue] (1) circle (0.5pt);
				\filldraw[white, draw=blue] (om) circle (0.5pt);
				\filldraw[white, draw=blue] (bom) circle (0.5pt);
				
			\end{tikzpicture}
			\caption{The geometric representation of the roots of \(\tilde{K}_n\).}
			\label{fig}
		\end{figure}
	\end{proof}
	
	\section{Various irreducibility-related results}
	
	Theorem \ref{maintrick} already gives the location of some of the roots of \(K_{a,n}\) in the general case. Computer calculations suggest that other roots do not form a smooth curve. However, for every individual \(K_{a,n}\), their location is similar to the case \(a=-1\), namely, some circle-like curves symmetric about the point \(\frac12\). We wish to explicitly find such a curve (which will depend on \(a\) and \(n\)). It most probably will not be as simple as in the case \(a=-1\) for the following reason: if we find a simple curve, it will generally yield a simple intersection with the real line. If we check that the intersection is not a root or not rational, we will get an elementary solution to FLT. 
	
	We anticipate that a square-freeness analysis similar to the case \(a=-1\) will be necessary.
	
	\begin{theorem}
		\(K_{a,n}\) is square-free for \(a\in\mathbb{Q}\setminus\{\pm1\}\).
	\end{theorem}
	We will denote \(A(x,y)\) the homogenization of the univariate rational polynomial \(A\) (i.e. \(A(x,y) = y^{\deg A}A\left( \frac{x}{y}\right) \)) and the \(d\)-th cyclotomic polynomial as \(\Phi_d\).
	\begin{lemma}
		\(\Phi_d(x,1-x)\) is irreducible in \(\mathbb{Q}[x]\). Moreover, \(\deg \Phi_d(x,1-x) = \varphi(d)\) whenever \(d\neq 2\).
	\end{lemma}
	\begin{proof}
		The statement is trivial for \(d=2\). Otherwise, consider any decomposition
		\begin{equation}\label{eq:2decompose}
			(1-x)^{\varphi(d)}\Phi_d\left(\frac{x}{1-x} \right) = A(x)B(x),\quad \deg A+\deg B =  \deg \Phi_d(x,1-x).
		\end{equation} 
		By making a change of variables \(x = \frac{t}{1+t}\), we will have
		\begin{multline*}
			\dfrac{1}{(1+t)^{\varphi(d)}}\Phi_d(t) = A\left(\frac{t}{1+t}\right)B\left(\frac{t}{1+t}\right) \implies\\ \Phi_d(t) = (1+t)^{\varphi(d)-\deg A - \deg B}A(t,1+t)B(t,1+t).
		\end{multline*}
		
		It is trivial that \(\deg \Phi_d(x,1-x)\leq \deg \Phi_d(x) = \varphi(d)\), so \eqref{eq:2decompose} implies that the power of \(1+t\) is non-negative. However, since cyclotomic polynomials are irreducible and \(\Phi_d(t)\neq t+1\), \(\gcd(\Phi_d(t), 1+t)=1\) . Hence, we have \(\deg A + \deg B = \deg \Phi_d(x,1-x) = \varphi(d)\) (the second part of the lemma) and \begin{equation}\label{eq:cycdecompose}
			\Phi_d(t) = A(t,1+t)B(t,1+t).
		\end{equation}
		Now we have \(\deg A(t,1+t) + \deg B(t,1+t) = \varphi(d) = \deg A + \deg B\). Since \(\deg A(t,1+t) \leq \deg A \) and \(\deg B(t,1+t) \leq \deg B \), we get \(\deg A(t,1+t) = \deg A\) and \(\deg B(t,1+t) = \deg B\). Note that \eqref{eq:cycdecompose} is a decomposition of a cyclotomic polynomial, so we must have \(\deg A = 0\) or \(\deg B = 0\), which means the decomposition \eqref{eq:2decompose} is trivial.
	\end{proof}
	\begin{proof}[Proof of main result]
		Let \[G=\gcd(K_{a,n} K'_{a,n})=\gcd(x^n+(1-x)^n+a^n, x^{n-1}-(1-x)^{n-1}).\]
		Similar to \eqref{eq:modtrick}, we can make a "modular arithmetic" argument and infer that \(G\mid x^n + x^{n-1}(1-x) + a^n = x^{n-1}+a^n\). The roots of \(G\) lie on the circle \(|x|=a^{\frac{n}{n-1}}\). On the other hand, they lie on the curve \(|x|=|1-x|\) as \(G\mid x^{n-1} - (1-x)^{n-1}\). If \(|a|^{\frac{n}{n-1}}<\frac12\), the intersection of the line and the circle is empty, so \(G=1\). We cannot have \(|a|^{\frac{n}{n-1}}=\frac12\) as \(|a|\) is rational and \(\left( \frac12 \right)^{n-1} \) is not the \(n\)th power of a rational. Therefore, we consider the case when the intersection contains two points. Denote those points \(B\) and \(\bar{B}\).
		
		Note that \(G\) has real coefficients. Moreover, \(G\mid x^{n-1}+a^n\) and the latter is square-free. Therefore, there are two possibilities: either \(G=1\) or \(G=(x-B)(x-\bar{B})=x^2-2x\Re B + |B|^2=x^2-x+|a|^{\frac{2n}{n-1}}\).
		
		Now consider the polynomial \(x^{n-1} - (1-x)^{n-1}\). Its canonical factorization is \[
		x^{n-1} - (1-x)^{n-1} = \prod_{d\mid n-1} \Phi_d(x,1-x).
		\]
		Assume \(G\neq 1\), then \(G\) is irreducible in \(\mathbb{R}[x]\). It must be equivalent (obtained by multiplication by a nonzero constant) to some polynomial among \(\Phi_d(x,1-x)\). Since \(\deg G = 2\), \(d\in\{3,4,6\}\), meaning \(x^2-x+|a|^{\frac{2n}{n-1}}\) is equivalent to one of 
		\[
		\Phi_3(x,1-x) = x^2-x+1, \Phi_4(x,1-x) = 2x^2-2x +1, \Phi_6(x,1-x) = 3x^2-3x+1.
		\]
		However, \(|a|^{\frac{2n}{n-1}} \neq 1\) as \(|a|\neq 1\), and \(|a|^{\frac{2n}{n-1}} \in \left\{\frac12,\frac13\right\}\) is impossible as \(|a|\) is rational.
	\end{proof}
	
	In the case \(a=1\), it is only meaningful to consider odd \(n\) as otherwise \(K_{1,n} = K_n\).
	
	\begin{theorem}
		\(K_{1,n}\) is squarefree for odd \(n\).
	\end{theorem}
	\begin{proof}
		Everything in the proof of the previous result (except the last line) is applicable to this case as well, so we have three candidates for \(G\) if \(G\neq 1\): \(G_1:=x^2-x+1, G_2:=2x^2-2x+1, G_3:=3x^2-3x+1\). These three polynomials have content \(1\), and \(K_{1,n}\in\mathbb{Z}[x]\), so \(G(x)\mid K_{1,n}(x)\) in \(\mathbb{Z}[x]\) by Gauss's Lemma. Plugging in any value for \(x\) must yield a correct divisibility relation in \(\mathbb{Z}\). However, \(K_{1,n}(2) = 2^n\), but \(G_1(2) = 3, G_2(2) = 5, G_3(2) = 7\), none of which divides \(2^n\). This contradiction proves \(G=1\).
	\end{proof}
	
	\begin{theorem}
		If \(a = \pm\frac12\), \(K_{a,n}\) is irreducible.
	\end{theorem}
	\begin{proof}
		Observe that \(L:= 2^nK_{a,n}\left( \frac12x\right) = x^n + (2-x)^n + (\pm1)^n\), and the irreducibility of \(L\) over \(\mathbb{Q}\) is equivalent to that of \(K_{a,n}\). Note that \(L\) has integer coefficients, so we will prove that it's irreducible over \(\ZZ\).
		
		The leading coefficient of \(L\) is either \(2\) (if \(n\) is even) or \(2n\) (if \(n\) is odd). In either case, it is divisible by \(2\) and not \(4\). Now assume \(L(x)=A(x)B(x)\). Without loss of generality, \(A\) has an odd leading coefficient. Considering this equation in the finite field \(\mathbb{F}_2\) yields \(A(x)B(x) = 1\), so they are both constant polynomials. However, the degree of \(A\) is preserved when passing to \(\mathbb{F}_2\) due to its odd leading coefficient, so \(A\) is constant. This means that \(L\) is irreducible.
	\end{proof}
	
	\begin{theorem}
		Let \(n\) be even, square-free, or square of a prime. For any irreducible factor \(P\in\ZZ[x]\) of \(\tilde K_n\) and every root \(r\) of \(P\), the symmetric copies \(1-r, \frac1r, \ldots\) (as depicted in Figure \ref{fig}) are also roots of \(P\). As a consequence, \(6\mid\deg P\).
	\end{theorem}
	\begin{proof}
		It is enough to prove the theorem for the case when \(P\) contains a root on the right arc. The problem will then follow by the following reasoning: any irreducible factor \(Q\) of \(\tilde K_n\) has a root \(r'\) which has a symmetric copy \(r''\) on the right arc. The minimal polynomial \(R\) of \(r''\) has all the symmetric copies of \(r''\), including \(r'\). Since \(\tilde K_n\) is square-free, \(Q\) has to be \(R\), so it satisfies the required property.
		
		Now let \(P\) be an irreducible factor of \(\tilde K_n\) which has a root \(r\) on the right arc. Since \(P\) has real coefficients, \(\bar r = \frac1r\) is also a root of \(P\). Polynomials \(P\) and \(P^*\) have a common root, so they are not relatively prime. Since \(P\) is irreducible, we get \(P\mid P^*\). But \(P\) and \(P^*\) have the same degree, so \(P^* = cP\) for some constant \(c\). Note that \(P^*(1) = 1^{\deg P}P(\frac11) = P(1)\implies c=1\).
		
		Next, observe that if \(P\) has a root on the line, similar reasoning would yield \(P(x)=P(1-x)\). \(P\) has two symmetries generating \(D_3\), so it satisfies the required condition. Similar reasoning applies if \(P\) has a root on the left arc (this time with the symmetry \(x\mapsto \frac{x}{x-1}\)). Thus, \(P\) has roots only on the right arc. We now proceed to showing that this case is impossible.
		
		Note that \(P(x), P(1-x), P(1-x)^*\) are all irreducible polynomials with disjoint set of roots, so \(P(x)P(1-x)P(1-x)^*\mid \tilde K_n\). Plug in \(x=0\) to get \(P(0)P(1)l \mid \tilde K_n(0)\), where \(l\) is the leading coefficient of \(P(1-x)\). Since \(P\) has no real roots, its degree is even, so \(l\) is also the leading coefficient of \(P(x)\). Considering the fact that \(P=P^*\), \(l=P(0) \implies P(0)^2P(1)\mid K_n(0)\). If \(n\) is even, \(K_n(0) = 2\), otherwise \(K_n(0) = n\). If \(n\) is even or squarefree, we get \(P(0)=\pm 1\). However, \[P(1) = l\prod_{\rho \text{ roots of } P} (1-\rho) \implies |P(1)| = |P(0)|\prod_{\rho } |1-\rho| < |P(0)|\prod_{\rho } 1 = 1,\]
		which is impossible.
		
		In the case \(n = p^2\), note that \(K_n\equiv 0\) in \(\mathbb{F}_p\), so all the coefficients of \(K_n\) are divisible by \(p\). Since \(P(x),P(1-x),P(1-x)^*\) are all primitive polynomials, we actually have \(P(x)P(1-x)P(1-x)^*\mid \frac1p \tilde K_n\). The value of the RHS at \(0\) is \(p\), a square-free number, so the same reasoning as before applies here.
	\end{proof}

	\section*{Acknowledgments}
	
	I express my sincere gratitude to Mihran Papikian whose comments helped me enormously as a novice researcher.


\begin{thebibliography}{0}
		
		\bibitem{inversion} H. S. M. Coxeter, {\it Introduction to geometry}, 2nd edition (Wiley, 1969).
		
		\bibitem{Edwards1996} Harold M. Edwards, {\it Fermat's last theorem: A genetic introduction to algebraic number theory} (Springer, 1996).
		
	\end{thebibliography}
\end{document}